\documentclass[12pt, reqno]{amsart}
 \usepackage{amsmath, amsthm, amscd, amsfonts, amssymb, graphicx, color, float}
\usepackage[bookmarksnumbered, colorlinks, plainpages]{hyperref}
\usepackage[utf8]{inputenc}
\input{mathrsfs.sty}
\hypersetup{colorlinks=true,linkcolor=red, anchorcolor=green, citecolor=cyan, urlcolor=red, filecolor=magenta, pdftoolbar=true}
\usepackage{soul,cancel}
\usepackage{tikzsymbols}
\usepackage{textcomp}
\newtheorem{theorem}{Theorem}[section]

\newtheorem{proposition}[theorem]{Proposition}
\newtheorem{corollary}[theorem]{Corollary}

\theoremstyle{definition}
\newtheorem{question}{Question}[]
\newtheorem{definition}[theorem]{Definition}
\newtheorem{example}[theorem]{Example}

\newtheorem*{JacConj}{Jacobian Conjecture}
\theoremstyle{remark}
\newtheorem{remark}[theorem]{Remark}

\numberwithin{equation}{section}

\newcommand{\R}{\mathbb{R}} 
\newcommand{\N}{\mathbb{N}} 

\newcommand{\C}{\mathbb{C}}

\newcommand{\K}{\mathbb{K}}

\begin{document}

\title[\tiny Counterexamples to the Jacobian Conjecture from the viewpoint of lineability]{Counterexamples to the Jacobian Conjecture from the viewpoint of lineability}


\author[Bernal]{L.~Bernal-Gonz\'alez}
\address[L.~Bernal-Gonz\'alez]{\mbox{}\newline \indent Departamento de An\'{a}lisis Matem\'{a}tico \newline
\indent Facultad de Matem\'{a}ticas \newline \indent
	Instituto de Matem\'aticas Antonio de Castro Brzezicki (IMUS)\newline
\indent Universidad de Sevilla \newline
\indent Avenida Reina Mercedes s/n, 41012-Sevilla (Spain).}
\email{lbernal@us.es}

\author[Calderón]{M.C.~Calder\'on-Moreno}
\address[M.C.~Calder\'on-Moreno]{\mbox{}\newline \indent Departamento de An\'{a}lisis Matem\'{a}tico \newline
\indent Facultad de Matem\'{a}ticas \newline \indent
	Instituto de Matem\'aticas Antonio de Castro Brzezicki (IMUS) \newline
\indent Universidad de Sevilla \newline
\indent Avenida Reina Mercedes s/n, 41012-Sevilla (Spain).}
\email{mccm@us.es}


\author[Prado]{J.A.~Prado-Bassas}
\address[J.A.~Prado-Bassas]{\mbox{}\newline \indent Departamento de An\'{a}lisis Matem\'{a}tico \newline
\indent Facultad de Matem\'{a}ticas \newline \indent
	Instituto de Matem\'aticas Antonio de Castro Brzezicki (IMUS)\newline
\indent Universidad de Sevilla \newline
\indent Avenida Reina Mercedes s/n, 41012-Sevilla (Spain).}
\email{bassas@us.es}





\begin{abstract}
Algebraic genericity in the family of polynomial counterexamples to the recently solved Jacobian Conjecture for $N\ge3$ is investigated. The notion of affine-coneability is introduced and it is shown that this family is affine-coneable and 1-lineable, but not 2-lineable. Semigroupability is further studied with respect to coordinatewise multiplication and composition.
\end{abstract}


\maketitle


\section{Introduction}

\quad The Jacobian Conjecture is one of the classical problems in Algebraic Geometry. It is usually attributed to Keller \cite{Keller}, who formulated the general version in 1939, although the two-dimensional case had already appeared in a paper of Kraus in 1884 (see \cite{Kraus} and \cite{LORD}). The conjecture is usually stated over an arbitrary field ${K}$ of characteristic zero. In this note, however, we will mainly work on $\C$, where it can be stated as follows.

\begin{JacConj}
Let $N\ge2$ and $\vec F:\C^N \to \C^N$ be a polynomial map in $N$ complex variables, such that $\det(J\vec F)$ is a non-zero constant.
Then $\vec F$ is invertible (with polynomial inverse).
\end{JacConj}

\vskip -2pt

Here, by $J\vec F$ we mean the classical Jacobian matrix, that is, if $\vec{F} = (F_1,\dots,F_N)$, then $J\vec F=\dfrac{\partial(F_1,\dots,F_N)}{\partial (z_1,\dots, z_N)}$.

\vskip 3pt

As one can see, the formulation of the Jacobian Conjecture is very simple. However, it has had a notable influence on affine Algebraic Geometry, so as to motivate many different approaches, reductions and equivalent formulations over the years. Its importance was highlighted by Smale, who included it as {\em Problem 16} in his list of mathematical problems for the twenty-first century \cite{Smale}.

\vskip 2.5pt

Very recently, on July 20, 2026, Alp\"oge \cite{Alpoge} announced an explicit counterexample in dimension 3,
obtained with the assistance of an AI (Claude Fable 5). This example easily extends to every dimension
$N > 3$ (just define the $j$-component of $F$ as the $j$th coordinate for $j \ge 4$).
The counterexample is very simple and it can be verified by a direct computation. Specifically, consider the vector function
\[
\vec F_0 = (P,Q,R) : \mathbb C^3 \longrightarrow \mathbb C^3
\]
given by
\begin{eqnarray*}
P(x,y,z) &=&
(1+xy)^3z+y^2(1+xy)(4+3xy),\\
Q(x,y,z)&=&
y+3x(1+xy)^2z+3xy^2(4+3xy),\\
R(x,y,z)&=&2x-3x^2y-x^3z.
\end{eqnarray*}
Then
$\det (J\vec F_0)=-2$ and since
$$
\vec F_0\left(0,0,-\frac14\right)
=
\vec F_0\left(1,-\frac32,\frac{13}{2}\right)
=
\vec F_0\left(-1,\frac32,\frac{13}{2}\right)
=
\left(-\frac14,0,0\right),
$$
$\vec F_0$ is not injective.

\vskip 3pt

From now on, let us denote by $\mathcal{CJC}_N$ the family of all counterexamples to the Jacobian Conjecture, that is, the family of non-injective polynomial maps with non-zero constant determinant of its Jacobian.

\vskip 3pt

If $N=1$ then it is immediate that the conjecture is true and, in fact, the solution is the set of polynomials with degree $1$. The original Jacobian Conjecture over $\C^2$, posed by Kraus in 1884, remains open (a partial answer to it was given in 1985 by Appelgate and Onishi \cite{AppelgateO}). Hence, $\mathcal{CJC}_1=\varnothing$, $\mathcal{CJC}_N\ne\varnothing$ for every $N\ge3$, whereas $\mathcal{CJC}_2$ may or may not be empty.

\vskip 3pt

The aim of this paper is to study how large can $\mathcal{CJC}_N$ be for $N\ge3$, in particular we are interested in its algebraic genericity. To this end, we shall use a number of techniques and notions taken from the modern Theory of Lineability, whose purpose is to detect large algebraic structures inside nonlinear sets. We refer the interested reader to the monograph \cite{aronBPS} and the seminal survey paper \cite{bams} for a broader introduction to this topic.

\vskip 3pt

The basic notions we will use are the following. Let $X$ be a vector space, $A\subset X$, and $\alpha$ a cardinal number. We say that $A$ is {\em $\alpha$-lineable} if there exists a vector subspace $V\subset X$ with $\dim(V) = \alpha$ and $V \setminus \{0\} \subset A$. The supremum of all cardinals $\alpha$ for which $A$ is $\alpha$-lineable is usually called the {\em degree of lineability} of $A$.

\vskip 3pt

The Theory of Lineability has traditionally focused on proving the existence of large algebraic structures. In fact, many classical examples turn out to be maximal-lineable
(i.e.~$\dim (X)$-lineable) or even larger structures (such as algebras). By contrast, examples of nonlinear sets with very small degree of lineability are rather unusual. In this paper we show that the family $\mathcal{CJC}_N$ exhibits precisely this behavior. More specifically, we prove that its degree of lineability is equal to 1. At the same time, we shall see that such an extreme form of non-lineability coexists with other algebraic structures.

\vskip 3pt

The lack of large vector spaces does not necessarily imply the absence of
significant algebraic structures. In many situations where we have low dimension lineability or even its total failure, one may still find other related structures. Following this idea, several authors have introduced some new weaker notions in the Theory of Lineability.

\vskip 3pt

If $X$ admits a semigroup structure, we shall also consider semigroupability, introduced by the authors in \cite{BCP_Peano}. Given a cardinal number $\alpha$, we say that a subset $A\subset X$ is {\em $\alpha$-semigroupable} if there exists an $\alpha$-generated semigroup $G\subset A$, that is, a semigroup $G$ generated by a set of cardinality $\alpha$, that cannot be generated by any set of strictly smaller cardinality.

\vskip 3pt

In 2008 \cite{Cone}, Aizpuru et al.~introduced the concept of coneability.
Given a vector space $X$, a subset $A\subset X$, and a cardinal number $\alpha$, we say that $A$ is {\em $\alpha$-coneable} if there exists a positive 
cone generated by an $\alpha$-linearly independent family. Unlike lineability, coneability does not require the presence of a vector space inside the set; it only retains the scalar multiplication structure.

\vskip 3pt

If we now allow convex combinations, we arrive to the notion of convex-lineabilty introduced, among others, by the first author in 2019 \cite{ConvexLin}. Let $A$ be a subset of a vector space $X$, and $\alpha$ a cardinal number. Then $A$ is said to be {\em $\alpha$-convex lineable} if there exists an infinite linearly independent subset $S\subset X$ with $\operatorname{card}(S)=\alpha$, such that its convex hull $\operatorname{Conv}(S)$ is contained in $A$.

\vskip 3pt

Finally, by joining the last two notions, the first and the second authors (jointly with other researchers) introduced in 2024 \cite{BernalConeable} the next concept. Assume that $X$ is a vector space over $\R$, $\alpha$ is a cardinal number,
and $A \subset X$. Then $A$ is said to be {\em $\alpha$-convex-coneable} if there exists a linearly independent set $S \subset X$ with $\operatorname{card}(S) = \alpha$ such that $\R\cdot\operatorname{Conv}(S) \subset A\cup\{0\}$.

\vskip 3pt

This paper is organized as follows. The degree of lineability of our family of counterexamples to the Jacobian Conjecture is established in section 2. In section 3, we introduce and study the notion of affine-coneability. This property, as well as a number of related ones together with semigroupability, is analyzed in section 4 for $\mathcal{CJC}_N$. The paper concludes with several remarks and open questions about extensions to these results to more general fields and new lineability properties.

\section{Low dimensional lineability of counterexamples}

\quad For $N\ge3$, let $\vec F_0\in\mathcal{CJC}_N$. It is obvious that for every $\lambda\in\C^*:=\C\setminus\{0\}$, the vector field $\lambda\vec F_0$ also belongs to $\mathcal{CJC}_N$. In the lineability setting, we say that this family is 1-lineable.

\vskip 3pt

The next result states that the family of maps with non zero constant determinant of its Jacobian can never be 2-lineable.

\begin{theorem}\label{th:2lin}
  For every $N\in\N$, the family
  \[ \mathcal S_N= \left\{ \vec F:\mathbb C^N\to\mathbb C^N: \, \det(J\vec F)\in\mathbb C^* \right\}. \] is not 2-lineable.
\end{theorem}

\begin{proof}
  By the way of contradiction, let assume the existence of a linear subspace $W$ of dimension at most 2 such that $W\setminus\{0\}\subset\mathcal{S}_N$.
  Let $\vec F,\vec G\in W$ be linearly independent vectors and define $\vec H_w:=\vec F+w\cdot\vec G\in W\setminus\{0\}\subset \mathcal{S}_N$ for any $w\in\C$.

\vskip 3pt

By hypothesis, $\vec H_w\ne0$, so $\det(J\vec H_w) =: c(w)\in\C^*$. But, on the other hand, $c(w)=\det(J\vec H_w)=\det(J\vec F+w\cdot J\vec G)$, and this is a polynomial on $w$ of degree $N$; more precisely, its leading coefficient is $\det(J\vec G)$ which is a non-zero constant (in fact, since $\det(J\vec H_w)$ only depends on $w$, every coefficient of the polynomial must be constant).

\vskip 3pt

Now a simple application of the Fundamental Theorem of Algebra guarantees the existence of $w_0\in\C^*$  such that $c(w_0)=0$ (in fact, since $\vec F\in\mathcal{S}_N$, $w_0\ne 0$), and $\vec H_{w_0}\notin\mathcal{S}_N$.
\end{proof}

Since for $N\ge3$ we have that $\mathcal{CJC}_N\subset \mathcal{S}_N$, we immediately get the next corollary.

\begin{corollary}
The family \,$\mathcal{CJC}_N$ $(N\ge3)$ has degree of lineability exactly $1$.
\end{corollary}

Moreover, since Theorem \ref{th:2lin} is true for any positive integer, one can also say that the family $\mathcal{CJC}_2$ is either empty or has degree of lineability $1$.

\section{Affine coneability}

Motivated by the weaker notions of lineability discussed above, we now introduce the concept of affine-coneability, which seems to be new. Roughly speaking, affine-coneability allows almost but not all possible linear combinations on the vector space structure.

\begin{definition} \label{def:affcone}
Let $X$ be a vector space over $\K$, where $\mathbb K=\mathbb R$ or $\mathbb C$, let $A\subset X$, and let $\alpha$ be a cardinal number.
Given a set $S\subset X$, we define its \emph{affine cone} by
$$\operatorname{ACone}(S) := \K^*\cdot \operatorname{Aff}(S),$$
where $\K^*=\mathbb K\setminus\{0\}$ and
$$\operatorname{Aff}(S) = \left\{ \sum_{j=1}^M \lambda_j x_j:\, M\in\N,\, x_j\in S,\, \lambda_j\in\K,\, \sum_{j=1}^M\lambda_j=1 \right\}.$$
Equivalently,
$$\operatorname{ACone}(S) = \left\{ \sum_{j=1}^M \lambda_j x_j:\, M\in\N,\, x_j\in S,\, \lambda_j\in\K,\, \sum_{j=1}^M\lambda_j\ne0 \right\}.$$
We say that $A$ is $\alpha$-{\em affine-coneable} if there exists a linearly independent set $S\subset A$, with $\operatorname{card}(S)=\alpha$, such that
$\operatorname{ACone}(S)\subset A$.
\end{definition}

As usual, if the cardinal is omitted, we understand \(\alpha=\aleph_0\). The following chain of implications follows directly from the definitions:
\[
\alpha\text{-lineability} \Longrightarrow \alpha\text{-affine-coneability} \Longrightarrow \alpha\text{-convex-coneability} \Longrightarrow \alpha\text{-convex-lineability}.
\]
The reverse implications, except the first one, are known not to hold in general. In the present paper, we shall see that the first reverse implication also fails, but it is not difficult to present an easier example.

\begin{example}
Let $N\in\N$ ($N\ge2$), $X=\R^N$, $A:=\{(x_1,\dots,x_N)\in\R^N:\, x_1+\dots+x_N\ne0\}$ and $S:=\{\vec e_1,\dots,\vec e_N\}$, where $\vec e_j$ is the $j$-th canonical vector. It is clear that $\operatorname{span}(S) = \R^N$ and that $\operatorname{ACone}(S)=A$. So $A$ is $N$-affine-coneable. However, $A$ is not even 2-lineable. Indeed, since $A=\R^N\setminus H$, where $H$ is a hyperplane of $\R^N$, every two-dimensional subspace $V$ of $\R^N$ has nontrivial intersection with $H$, so $V\setminus\{0\}\not\subset A$.
\end{example}

\begin{remark}
If $S$ is linearly independent and $V=\operatorname{span}(S)$, then $\operatorname{ACone}(S)$ can be viewed as the complement of a hyperplane in $V$. Indeed, every $x\in V$ has a unique representation
$$x=\sum_{j=1}^M \lambda_j x_j \qquad (x_j\in S),$$
and the map
$T:V\to \K$ given by $$T\left(\sum_{j=1}^M\lambda_jx_j\right)=\sum_{j=1}^M\lambda_j $$
is a nonzero linear functional. Hence
$$\operatorname{ACone}(S)=V\setminus\ker(T).$$
Thus, affine-coneability may be viewed as {\em a hyperplane punctured version} of li\-nea\-bi\-li\-ty: instead of containing a vector space except for the origin, the set contains a vector space except for a whole hyperplane.
\end{remark}

In fact, this is an equivalent definition of affine-coneability:

\begin{proposition}
  Let $X$ be a vector space over $\K$, let $A\subset X$, and let $\alpha$ be a cardinal number. The set $A$ is $\alpha$-affine-coneable if and only if there exist a vector subspace $V\subset X$, with $\dim(V)=\alpha$, and a hyperplane $H\subset V$ such that $V\setminus H\subset A$.
\end{proposition}

\begin{proof}
The only if part is just the above remark. Let us proof the converse.

\vskip 3pt

Assume that there exist a vector subspace $V\subset X$ with $\dim(V)=\alpha$ and a hyperplane $H\subset V$ such that $V\setminus H\subset A$. Since $H$ is a hyperplane of $V$, there is $v\in V\setminus H$ such that $\operatorname{span}\{v\}\oplus H=V$. Let $S_0$ be a basis of $H$ and take $S:=\{v\}\cup\{v+h:\, h\in S_0\}$. By construction, we have $S\subset V\setminus H$ and $\operatorname{card}(S)=\dim(V)=\alpha$. Moreover, $S$ is linear independent. Indeed, if for $M\in\N$, $h_1,\dots,h_M\in H$ and $\lambda_0,\lambda_1,\dots,\lambda_M\in\K$ we have
$$\lambda_0 v+\sum_{j=1}^M \lambda_j(v+h_j)=0,$$
then
$$ \left(\sum_{j=0}^M\lambda_j\right)v + \sum_{j=1}^M\lambda_j h_j=0.$$
But $V=\operatorname{span}\{v\}\oplus H$, so $$\sum_{j=0}^M\lambda_j=0 \quad\text{and}\quad \sum_{j=1}^M\lambda_jh_j=0.$$ Now, recall that $S_0$ is a basis of $H$; this implies that $\lambda_j=0$ ($1\le j\le M$), and then also $\lambda_0 = 0$.

\vskip 3pt

Finally, let us show that $\operatorname{ACone}(S)\subset V\setminus H$. As seen before, every finite linear combination of elements of $S$ is of the form $$\left(\sum_{j=0}^M\lambda_j\right)v+\sum_{j=1}^M\lambda_j h_j.$$ By the direct decomposition $V = \operatorname{span}\{v\}\oplus H$, this vector belongs to $H$ if and only if the coefficient of $v$ is zero, that is, $\sum_{j=0}^M\lambda_j=0$. This concludes the proof.
\end{proof}

As the next example shows, \(\alpha\)-affine-coneability can be seen as a property lying immediately below \(\alpha\)-lineability in the hierarchy of algebraic genericity notions.

\begin{example}\label{ex:constr}
  Given a cardinal number $\alpha\ge2$, there is a vector space $X$ (with dimension at least $\alpha$) and a subset $A\subset X$ such that
  \begin{itemize}
    \item $A$ is $\alpha$-affine-coneable.
    \item $A$ is not $\alpha$-lineable.
    \item $A$ is $\beta$-lineable for all cardinal $\beta<\alpha$.
  \end{itemize}

\vskip 3pt

Let $X$ be a vector space with dimension large enough (for example, $\dim(X)\ge\sum_{\beta\le \alpha}\beta\ge\alpha$), so that for each cardinal $\beta\le\alpha$, there exist a vector subspace $V_\beta$ of $X$ with $\dim(V_\beta)=\beta$, and such that $V_\beta\cap V_{\beta'}=\{0\}$ for any pair $\beta\ne\beta'$.

\vskip 3pt

Let $H\subset V_\alpha$ be a hyperplane and let $A := (V_\alpha\setminus H)\cup \bigcup_{\beta<\alpha} \left( V_\beta \setminus \{0\} \right)$. Observe that every element of $A$ belongs to a unique $V_\beta$ with $\beta\le\alpha$. Furthermore, no element of $A$ can be decomposed as a sum of two (or more) nonzero elements of from distinct $V_\beta$'s.

\vskip 3pt

It is trivial that $A$ is $\alpha$-affine-coneable and $\beta$-lineable for every $\beta<\alpha$. Let us show that $A$ cannot be $\alpha$-lineable. For that purpose, let $W\subset X$ be a vector subspace such that $W\setminus\{0\}\subset A$. Then there is a unique $\beta_0\le\alpha$ such that $W \subset V_{\beta_0}$.

\vskip 3pt

Indeed, let $w\in W\setminus \{0\}$. Since $w\in A$, there is a unique $\beta_0\le\alpha$ such that $w\in V_{\beta_0}$. Suppose that there is $w'\in W\setminus V_{\beta_0}$; since $w'\in A$, then again there is (a unique) $\beta'\le\alpha$ with $\beta'\ne\beta_0$ such that $w'\in V_{\beta'}$. In particular, $w$ and $w'$ are linearly independent and $w+w'\ne0$ is an element of $X$ with nonzero components in two different subspaces $V_\beta$, so $w+w'\notin A$. But since $w,w'\in W$, we have that $w+w'\in W\setminus\{0\}\subset A$, which is a contradiction.

\vskip 3pt

Finally, if $\beta_0<\alpha$, then $\dim(W)\le\dim(V_{\beta_0})=\beta_0<\alpha$. On the other hand, if $\beta_0=\alpha$, we have that $W\setminus\{0\}\subset V_\alpha\setminus H$, hence $W\cap H=\{0\}$. But $H$ is a hyperplane of $V_\alpha$, so $\dim(W)\le1<2\le\alpha$.
\end{example}

\section{Affine-coneability and further algebraic structures of counterexamples}

In section 2 we proved that, for $N\ge 3$, the family $\mathcal{CJC}_N$ is 1-lineable but not 2-lineable. In particular, it cannot contain any nontrivial linear structure of dimension greater than one. However, the notion of affine-coneability introduced in the previous section shows that this lack of lineability does not prevent the existence of large algebraic structures.

\vskip 3pt

In this section we prove that such structures occur in a pointwise way. In order to be more precise,
and inspired by the concept of pointwise lineability coined by Pellegrino and Raposo in \cite{pellegrinoraposo},
we introduce the following notion, which is trivially stronger than the ordinary $\alpha$-affine-coneability.

\begin{definition}
We say that a subset $A$ of a vector space $X$ is {\em pointwise \(\alpha\)-affine-coneable} if, for every $a\in A$, there exists a linearly independent set $S=S(a)\subset A$, with $\operatorname{card}(S)=\alpha$ and $a\in \operatorname{ACone}(S)\subset A$.
\end{definition}

Under this terminology, the main result of this section reads as follows.

\begin{theorem}
  The family $\mathcal{CJC}_N$ ($N\ge3$) is pointwise affine-coneable.
\end{theorem}

\begin{proof}
Let $\vec F_0\in\mathcal{CJC}_N$, and denote by $F_{0i}$ ($1\le i\le N$) its components.
Let $\{p_j\}_{j\in\N}$ be a basis of the vector space of all polynomials of $N-1$ complex variables, and denote $p_0 \equiv 0$.

\vskip 3pt

For every $j \ge 0$, let $\vec G_j : \C^N \longrightarrow \C^N$ be vector field given by
$$
\vec G_j(z_1,\dots,z_N) := (p_j(z_2,\dots,z_N),0,\cdots,0),
$$
and define $\vec H_j := \vec F_0 + \vec G_j(\vec F_0)$. Also, let $S := \{\vec H_j:\, j\ge0\}$. Observe that $\vec F_0=\vec H_0\in S$, just because $p_0 \equiv 0$.

\vskip 3pt

First, we claim that $S\subset \mathcal{CJC}_N$. Indeed, since $F_{0i}$ ($1\le i\le N$) and $p_j$ ($j\ge0$) are polynomials, it is immediate that the components of every $\vec G_j$ (and, thus, $\vec H_j$) are also polynomials. Moreover a simple computation leads us to
$\vec H_j = \vec\Phi_j\circ \vec F_0$, where $\vec\Phi_j:\C^N\to\C^N$ is the vector field given by
$$\vec \Phi_j(w_1,w_2,\dots,w_N):=(w_1+p_j(w_2,\dots,w_N),w_2,\dots,w_N).$$ But it is trivial that $\det(J \vec\Phi_j)=1$, so, using that the Jacobian of the composition is the product of the Jacobians, we get that
$\det(J \vec H_j) = \det(J(\vec F_0))\in\C^*$. It only remains to show that $\vec H_j$ is not injective, but this comes from the non injectivity of $\vec F_0$ and the fact that $\vec H_j = \vec \Phi_j\circ \vec F_0$.

\vskip 3pt

Let us show now that the family $S$ is linearly independent. For that purpose, let $M\in\N\cup\{0\}$ and $\lambda_0,\dots,\lambda_M\in\C$ such that $$\sum_{j=0}^M \lambda_j\vec H_j=\vec 0.$$
Taking into account that $\vec H_j=\vec F_0+\vec G_j$ and that $\vec G_0=\vec 0$, if we denote by $\Lambda:=\sum_{j=0}^M\lambda_j$, we get that
$$
\Lambda\cdot\vec F_0+ \left( \sum_{j=1}^M\lambda_j p_j(F_{02},\dots,F_{0N}),0,\dots,0\right) = \vec0.
$$

For $2\le j\le N$, we get that $\Lambda F_{0j}=0$, so $\Lambda = 0$ (if not, $F_{0j}=0$ for $2 \le j \le N$ and $J \vec F_0$ would have $N-2$ null rows, so zero determinant, which would exclude $\vec F_0$ from $\mathcal{CJC}_N$).
Turning now our attention to the first component, we obtain that
\begin{equation}\label{eq:1}\sum_{j=1}^M\lambda_jp_j(F_{02},\dots,F_{0N})= 0.\end{equation}
But by the Jacobian criterion for algebraic independence (see, for instance \cite{Lef}), the nonvanishing of $\det(J\vec F_0)$ implies that all its components are algebraically independent, so \eqref{eq:1} is equivalent to $$\sum_{j=1}^M\lambda_jp_j(z_2,\dots,z_N)=0.$$ Now, observe that $\{p_j\}_{j\ge1}$ is a basis of polynomials of $N-1$ variables, so $\lambda_j=0$ for $j=1,\dots,M$, and recalling that $\Lambda=0$, we also get $\lambda_0 = 0$.

\vskip 3pt

Since $\vec F_0\in S\subset \operatorname{ACone}(S)$, it only remains to show that $\operatorname{ACone}(S) \subset \mathcal{CJC}_N$. Let $M\in\N$, and let  $\lambda_0,\lambda_1\dots,\lambda_M\in\C$ with $\Lambda:=\sum_{j=1}^M\lambda_j\ne0$. Observe that
$$\sum_{j=0}^M\lambda_j \vec H_j =\Lambda \vec F_0+(p(F_{02},\dots,F_{0N}),0,\dots,0),$$ where $$p(z_2,\dots,z_N)=\sum_{j=1}^M\lambda_j p_j(z_2,\dots,z_N).$$
Since $\Lambda\ne0$, we have that $\Lambda\vec F_0\in\mathcal{CJC}_N$. Proceeding as in the proof of $S\subset\mathcal{CJC}_N$, we are done.
\end{proof}

Observe that, since the vector space of all polynomial maps over $\C^N$ has infinite countable dimension, the above result is sharp in the sense of cardinality.

\vskip 3pt

We conclude this section by turning our attention to further algebraic stability pro\-per\-ties of $\mathcal{CJC}_N$, namely its behavior under coordinatewise multiplication and composition.

\begin{proposition}
The family $\mathcal{CJC}_N$ is not 1-semigroupable under coordinatewise pro\-duct.
\end{proposition}

\begin{proof}
It suffices to show that if $\vec F_0\in\mathcal{CJC}_N$, then $\vec F_0^2\notin \mathcal{CJC}_N$.

\vskip 3pt

Let $\vec F_0\in\mathcal{CJC}_N$ and write $c:=\det (J\vec F_0)\in\C^*$. Since $\vec F_0^2=(F_{01}^2,\dots,F_{0N}^2)$, it is a simple computation that $J\vec F_0^2=\operatorname{diag}(2F_{01},\dots,2F_{0N})\cdot J\vec F_0$, whence
$$\det(J\vec F_0^2)=2^Nc\prod_{i=1}^N F_{0i}.$$
Now, since $\vec F_0\in\mathcal{CJC}_N$, its components are polynomials and none of them can be constant (otherwise, we get a null-row in its Jacobian), so the term $\displaystyle \prod_{i=1}^N F_{0i}$ is a non constant polynomial, which implies that $\vec F_0^2\notin\mathcal{CJC}_N$.
\end{proof}

Observe now that a semigroup generated by a set of cardinal $\alpha\ge\aleph_0$ has cardinal at most $\alpha$, since every element can be written as a finite product of elements of the generating set and there are only countably many such products. Consequently, every semigroup of cardinality $\alpha>\aleph_0$ requires a generating set of cardinality  $\alpha$. Since the semigroup itself is a generating set, it admits a generating set of cardinality exactly $\alpha$. So, in particular, every semigroup of cardinality $\mathfrak c$ is $\mathfrak c$-semigroupable.

\begin{proposition}
  The family \,$\mathcal{CJC}_N$ is a semigroup, under composition, of cardinal $\mathfrak{c}$. In particular it is $\mathfrak{c}$-semigroupable.
\end{proposition}

\begin{proof}
The family $\mathcal{CJC}_N$ is stable under composition. Indeed, composition of polynomials is again a polynomial; the Jacobian of a composition is the product of Jacobians, so the product of two non-zero constant determinants of Jacobians gives another non-zero constant determinant; lastly, $\vec F\circ\vec G$ is non-injective whenever $\vec G$ is non injective (and $\vec F$ a map). Therefore it forms a semigroup.

\vskip 2pt

Finally, given $\vec F\in\mathcal{CJC}_N$, the set $A := \{\vec F+(a,0,\dots,0):\, a\in\C\}$ satisfies that $A\subset \mathcal{CJC}_N$ and $\operatorname{card}(A) =\operatorname{card}(\C) = \mathfrak{c}$. Since the set of all polynomial maps $\C^N\to\C^N$ has cardinality $\mathfrak c$, it follows that $\mathcal{CJC}_N$ has cardinality exactly $\mathfrak{c}$.
\end{proof}

\section{Final remarks and open problems}

\noindent {\bf 1.} The interpretation of affine-coneability as a hyperplane-punctured form of lineability naturally leads to a broader family of punctured lineability notions. Namely, one may replace the removed hyperplane by a subspace of higher codimension. To be more precise, we introduce the following concept.

\begin{definition}
If $1\le \beta\le \alpha$ are cardinal numbers, we say that a set $A\subset X$ is {\em $\alpha\ominus\beta$-punctured-lineable} if there exist an $\alpha$-dimensional vector subspace $V\subset X$ and a subspace $Y\subset V$ with $\operatorname{codim}(Y)=\beta$ such that $V\setminus Y\subset A$.
\end{definition}

With this terminology, $\alpha$-affine-coneability is the same as $\alpha\ominus1$-punctured-lineability, while $\alpha\ominus\alpha$-punctured lineability is, in fact, $\alpha$-lineability. In addition, it is straightforward that $\alpha\ominus\beta_2$-punctured-lineability implies $\alpha\ominus\beta_1$-punctured-lineability, for cardinals $1\le \beta_1<\beta_2\le \alpha$.

\vskip 3pt

Moreover, a construction similar to that of Example~\ref{ex:constr} shows that, given cardinals $1\leq \beta_1<\beta_2\leq \alpha$, there exists a subset $A$ of a vector space $X$ which is $\alpha\ominus\beta_1$-punctured-lineable but not $\alpha\ominus\beta_2$-punctured-lineable.

\vskip 3pt

In view of the results obtained for the family of counterexamples to the Jacobian Conjecture, the following question arises naturally.

\begin{question} Let $N\geq 3$. To what extent does the family $\mathcal{CJC}_N$ satisfy higher punctured-lineability properties?

\noindent More precisely, is $\mathcal{CJC}_N$ (pointwise) $\aleph_0\ominus n$-punctured-lineable for some $n\in\N$ ($n>1$)?
\end{question}

\vskip 5pt 

\noindent {\bf 2.} The Jacobian Conjecture can be formulated over a general field \,$K$ \,with characteristic $0$. However, all these versions are known to be equivalent, so it is enough to establish (or refute) the conjecture over $\C$, and the counterexample given in \cite{Alpoge} implies that the family $\mathcal{CJC}_N(K)$ of counterexamples of the Jacobian Conjecture over $K^N$ is nonempty for every $N\ge3$ and every field $K$ of characteristic 0.

\vskip 3pt

Most of the results in this paper have been stated over $\C$, but they remain valid over arbitrary fields of characteristic \(0\), with exactly the same proofs. The only exception is Theorem~\ref{th:2lin}, whose proof relies on the existence of a root of a certain polynomial of degree $N$. Hence the same argument remains valid over every algebraically closed field of characteristic 0. Moreover, if $N$ is odd, it also holds over every real closed field, and therefore over $\R$. So we can state a more general result (see Theorem \ref{thm general} below).

\vskip 3pt

Prior to this, it should be emphasized that the equivalence of the Jacobian Conjecture over different fields of characteristic 0 concerns only the existence of counterexamples. It does not imply that the sets $\mathcal{CJC}_N(K)$ have similar algebraic properties. In particular, lineability questions may depend on the arithmetic properties of the ground field, as the next result illustrates.

\begin{theorem} \label{thm general}
Let $K$ be a field of characteristic 0. Let $\mathcal{CJC}_N(K)$ be the family of all counterexamples to the Jacobian Conjecture over $K^N$. We have:
\begin{itemize}
  \item For every $N\ge3$, the family $\mathcal{CJC}_N(K)$ is 1-lineable, $\operatorname{card}(K)$-semigroupable under composition, and pointwise affine-coneable, but it is not semigroupable under coordinatewise product.
  \item If \,$K$ is algebraically closed, then $\mathcal{CJC}_N(K)$ is not 2-lineable for every $N\ge3$.
  \item If \,$K$ is a real closed field, then $\mathcal{CJC}_N(K)$ is not 2-lineable for every odd number $N\ge3$.
\end{itemize}
\end{theorem}

We conclude the paper with the following natural question.

\begin{question}
Is $\mathcal{CJC}_N(\R)$ also not 2-lineable for every {\it even} $N\ge 4$?
\end{question}

\vskip 6pt

\noindent {\bf Acknowledgments.} The authors thank IMUS-Mar\'{\i}a de Maeztu grant CEX2024-001517-M - Apoyo a Unidades de Excelencia Mar\'{\i}a de Maeztu for supporting this research, funded by MICIU/AEI/10.13039/501100011033. The first and third author thank the Projects Programa operativo PPIT-FEDER Andaluc\'{\i}a 2021-2027 SOL2024-31596 and  the first and second author also thank the project SOL\-2024-31708 for partially supporting this research. The three authors also thank Miguel \'Angel Morales for bringing the existence of the counterexample to the Jacobian Conjecture to their attention (see \cite{Gauss}).

\vskip 6pt

\noindent {\bf Authors contribution.} All authors have contributed equally to this work.

\end{document}